\documentclass[10pt,headrule,fleqn,onecolumn]{article}

\usepackage{CJK} 
\usepackage{amssymb}
\usepackage{amsmath}
\usepackage{amsfonts}
\usepackage{amsopn}
\usepackage{amscd}
\usepackage{amsthm}
\usepackage{amsfonts,amssymb,color}
\usepackage{bm}
\usepackage{cite}
\usepackage{curves}
\usepackage{CJK}
\usepackage{delarray}
\usepackage{dsfont}
\usepackage{exscale}
\usepackage{framed}
\usepackage{hyperref}
\usepackage{indentfirst,latexsym,bm}
\usepackage{mathrsfs}
\usepackage{pifont}
\usepackage{xypic}

\usepackage{geometry}
\date{}
\begin{document}
\begin{CJK*}{GBK}{song}
\title{\bfseries\scshape All fractional $(g,f)$-factors in graphs
\thanks{Supported by the National Natural Science Foundation of China (Grant No. 11371009) and
the National Social Science Foundation of China (Grant No. 14AGL001).}
}
\author{\small Zhiren Sun$^{1}$, Sizhong Zhou$^{2}$\footnote{Corresponding
author. E-mail address: zsz\_cumt@163.com (S. Zhou)}\\
\small $1$. School of Mathematical Sciences, Nanjing Normal University, Nanjing 210046, P. R. China\\
\small $2$. School of Mathematics and Physics, Jiangsu University of Science and Technology,\\
\small Mengxi Road 2, Zhenjiang, Jiangsu 212003, P. R. China\\
}

\maketitle

\begin{abstract}
\noindent Let $G$ be a graph, and $g,f:V(G)\rightarrow N$ be two functions with $g(x)\leq f(x)$
for each vertex $x$ in $G$. We say that $G$ has all fractional $(g,f)$-factors if $G$ includes a
fractional $r$-factor for every $r:V(G)\rightarrow N$ such that $g(x)\leq r(x)\leq f(x)$ for
each vertex $x$ in $G$. Let $H$ be a subgraph of $G$. We say that $G$ admits all fractional
$(g,f)$-factors including $H$ if for every $r:V(G)\rightarrow N$ with
$g(x)\leq r(x)\leq f(x)$ for each vertex $x$ in $G$, $G$ includes a fractional $r$-factor $F_h$
with $h(e)=1$ for any $e\in E(H)$, then we say that $G$ admits all fractional $(g,f)$-factors including
$H$, where $h:E(G)\rightarrow [0,1]$ is the indicator function of $F_h$. In this
paper, we obtain a characterization for the existence of all fractional $(g,f)$-factors
including $H$ and pose a sufficient condition for a graph to have all fractional
$(g,f)$-factors including $H$.
\\
\begin{flushleft}

(2010) Mathematics Subject Classification: 05C70, 05C72
\end{flushleft}
\end{abstract}
\section{Introduction}
We consider finite undirected graphs which have neither multiple edges nor loops. Let $G$
be a graph. We denote its vertex set and edge set by $V(G)$ and $E(G)$, respectively. For each
$x\in V(G)$, the degree of $x$ in $G$ is defined as the number of edges which are adjacent to
$x$ and denoted by $d_G(x)$. For any $S\subseteq V(G)$, we use $G[S]$ to denote the subgraph
of $G$ induced by $S$, and use $G-S$ to denote the subgraph obtained from $G$ by deleting vertices
in $S$ together with the edges incident to vertices in $S$. A subset $S$ of $V(G)$ is said to be
independent if $N_G(S)\cap S=\emptyset$. Let $S$ and $T$ be two disjoint vertex subsets of $G$.
Then $e_G(S,T)$ denotes the number of edges joining $S$ to $T$.

Let $g,f:V(G)\rightarrow N$ be two functions with $g(x)\leq f(x)$ for each $x\in V(G)$. A spanning
subgraph $F$ of $G$ is called a $(g,f)$-factor if one has $g(x)\leq d_F(x)\leq f(x)$ for each
vertex $x$ in $G$. An $(f,f)$-factor is said to be an $f$-factor. If $G$ includes an $r$-factor
for every $r:V(G)\rightarrow N$ which satisfies $g(x)\leq r(x)\leq f(x)$ for each vertex $x$ in $G$
and $r(V(G))$ is even, then we say that $G$ admits all $(g,f)$-factors. Let $h:E(G)\rightarrow [0,1]$
be a function. For any $x\in V(G)$, we denote the set of edges incident with $x$ by $E(x)$. If
$g(x)\leq\sum_{e\in E(x)}h(e)\leq f(x)$ holds for each vertex $x$ in $G$, then we call graph $F_h$
with vertex set $V(G)$ and edge set $E_h$ a fractional $(g,f)$-factor of $G$ with indicator function
$h$, where $E_h=\{e:e\in E(G),h(e)>0\}$. A fractional $(f,f)$-factor is called a fractional $f$-factor.
If $G$ contains a fractional $r$-factor for every $r:V(G)\rightarrow N$ with $g(x)\leq r(x)\leq f(x)$
for each vertex $x$ in $G$, then we say that $G$ admits all fractional $(g,f)$-factors. If $g(x)\equiv a$, $f(x)\equiv b$ and $G$ admits all fractional $(g,f)$-factors, then we say that $G$ contains all fractional $[a,b]$-factors. Let $H$ be a subgraph of $G$. If for every $r:V(G)\rightarrow N$ such that
$g(x)\leq r(x)\leq f(x)$ for each vertex $x$ in $G$, $G$ includes a fractional $r$-factor $F_h$
with $h(e)=1$ for any $e\in E(H)$, then we say that $G$ admits all fractional $(g,f)$-factors including
$H$, where $h$ is the indicator function of $F_h$.  For any function $\varphi:V(G)\rightarrow N$, we
define $\varphi(S)=\sum_{x\in S}\varphi(x)$ and $\varphi(\emptyset)=0$. Especially, $d_G(S)=\sum_{x\in S}d_G(x)$.

Lu \cite{L} first introduced the definition of all fractional $(g,f)$-factors, and obtained a
necessary and sufficient condition for a graph to have all fractional $(g,f)$-factors, and posed a
sufficient condition for the existence of all fractional $[a,b]$-factors in graphs. Zhou and Sun
\cite{ZS} showed a neighborhood condition for a graph to have all fractional $[a,b]$-factors, which
is an extension of Lu's result \cite{L}. Zhou, Bian and Sun \cite{ZBS} obtained a binding number
condition for the existence of all fractional $[a,b]$-factors in graphs. The following results
on fractional $(g,f)$-factors and all all fractional $(g,f)$-factors are known.

\medskip

Anstee \cite{A} gave a necessary and sufficient condition for graphs to have fractional
$(g,f)$-factors. Liu and Zhang \cite{LZ} posed a new proof.

\medskip

\noindent{\textbf{Theorem 1}} (Anstee \cite{A}, Liu and Zhang \cite{LZ}). Let $G$ be a graph,
and $g,f:V(G)\rightarrow Z^{+}$ be two functions with $g(x)\leq f(x)$ for each
vertex $x$ in $G$. Then $G$ contains a fractional $(g,f)$-factor if and only if
$$
f(S)+d_{G-S}(T)-g(T)\geq0
$$
for any subset $S$ of $V(G)$, where $T=\{x:x\in V(G)-S, d_{G-S}(x)<g(x)\}$.

\medskip

The following theorem is equivalent to Theorem 1.

\medskip
\noindent{\textbf{Theorem 2}}. Let $G$ be a graph,
and $g,f:V(G)\rightarrow Z^{+}$ be two functions with $g(x)\leq f(x)$ for each
vertex $x$ in $G$. Then $G$ contains a fractional $(g,f)$-factor if and only if
$$
f(S)+d_{G-S}(T)-g(T)\geq0
$$
for all disjoint subsets $S$ and $T$ of $V(G)$.

\medskip

Lu \cite{L} showed a characterization of graphs having all fractional $(g,f)$-factors.

\medskip

\noindent{\textbf{Theorem 3}} (Lu \cite{L}). Let $G$ be a graph and $g,f:V(G)\rightarrow Z^{+}$
be two functions with $g(x)\leq f(x)$ for each vertex $x$ in $G$. Then $G$ admits all
fractional $(g,f)$-factors if and only if
$$
g(S)+d_{G-S}(T)-f(T)\geq0
$$
for any subset $S$ of $V(G)$, where $T=\{x:x\in V(G)-S, d_{G-S}(x)<f(x)\}$.

\medskip

In this paper, we study the existence of all fractional $(g,f)$-factors including any
given subgraph in graphs, and pose some new results which are shown in the following.

\medskip

\noindent{\textbf{Theorem 4}}. Let $G$ be a graph and $g,f:V(G)\rightarrow Z^{+}$ be two
functions such that $g(x)\leq f(x)$ for each vertex $x$ in $G$. Let $H$ be a subgraph of
$G$. Then $G$ has all fractional $(g,f)$-factors including $H$ if and only if
$$
g(S)+d_{G-S}(T)-f(T)\geq d_H(S)-e_H(S,T)
$$
for all disjoint subset $S$ and $T$ of $V(G)$.

\medskip

\noindent{\textbf{Theorem 5}}. Let $G$ be a graph, $H$ be a subgraph of $G$, and
$g,f:V(G)\rightarrow Z^{+}$ be two functions with $d_H(x)\leq g(x)\leq f(x)\leq d_G(x)$
for each vertex $x$ in $G$. If $(g(x)-d_H(x))d_G(y)\geq(d_G(x)-d_H(x))f(y)$ holds for any
$x,y\in V(G)$, then $G$ has all fractional $(g,f)$-factors including $H$.

\medskip

If $E(H)=\emptyset$ in Theorem 5, then we obtain the following corollary.

\medskip

\noindent{\textbf{Corollary 6}}. Let $G$ be a graph, and
$g,f:V(G)\rightarrow Z^{+}$ be two functions with $g(x)\leq f(x)\leq d_G(x)$
for each vertex $x$ in $G$. If $g(x)d_G(y)\geq d_G(x)f(y)$ holds for any
$x,y\in V(G)$, then $G$ contains all fractional $(g,f)$-factors.

\medskip

\section{The proof of Theorem 4}
\noindent{\it Proof of Theorem 4.} \ We first verify this sufficiency. Let $r:V(G)\rightarrow Z^{+}$
be an arbitrary integer-valued function such that $g(x)\leq r(x)\leq f(x)$ for each $x\in V(G)$. According
to the definition of all fractional $(g,f)$-factors including $H$, we need only to verify that $G$
admits a fractional $r$-factor including $H$, that is, we need only to verify that $G$ admits a
fractional $r'$-factor excluding $H$, where $r'(x)=d_G(x)-r(x)$. Let $G'=G-E(H)$. Thus, we need only
to prove that $G'$ admits a fractional $r'$-factor.

For any disjoint subsets $S$ and $T$ of $V(G)$,
$$
g(S)+d_{G-S}(T)-f(T)\geq d_H(S)-e_H(S,T),
$$
and so,
$$
g(T)+d_{G-T}(S)-f(S)-d_H(T)+e_H(S,T)\geq0.\eqno(1)
$$
It follows from (1) that
\begin{eqnarray*}
&&r'(S)+d_{G'-S}(T)-r'(T)=r'(S)+d_{G-S}(T)-r'(T)-d_H(T)+e_H(S,T)\\
&&=d_G(S)-r(S)+d_{G-S}(T)-d_G(T)+r(T)-d_H(T)+e_H(S,T)\\
&&\geq d_G(S)-f(S)+d_{G-S}(T)-d_G(T)+g(T)-d_H(T)+e_H(S,T)\\
&&=g(T)+d_{G-T}(S)-f(S)-d_H(T)+e_H(S,T)\geq0.
\end{eqnarray*}
In terms of Theorem 2, $G'$ admits a fractional $r'$-factor, that is, $G$ has all fractional
$(g,f)$-factors including $H$.

Now we verify the necessary. Conversely, we assume that there exist disjoint subsets $S$ and
$T$ of $V(G)$ such that
$$
g(S)+d_{G-S}(T)-f(T)<d_H(S)-e_H(S,T).
$$

Let $r(x)=g(x)$ for any $x\in S$ and $r(y)=f(y)$ for any $y\in V(G)\setminus S$. Thus, we have
\begin{eqnarray*}
0&>&g(S)+d_{G-S}(T)-f(T)-d_H(S)+e_H(S,T)\\
&=&r(S)+d_{G-S}(T)-r(T)-d_H(S)+e_H(S,T).
\end{eqnarray*}
Set $r'(x)=d_G(x)-r(x)$ and $G'=G-E(H)$. Thus,
\begin{eqnarray*}
0&>&r(S)+d_{G-S}(T)-r(T)-d_H(S)+e_H(S,T)\\
&=&d_G(S)-r'(S)+d_{G'-S}(T)+d_H(T)-e_H(S,T)-d_G(T)+r'(T)-d_H(S)+e_H(S,T)\\
&=&d_G'(S)+d_H(S)-r'(S)+d_{G'-S}(T)+d_H(T)-d_G'(T)-d_H(T)+r'(T)-d_H(S)\\
&=&r'(T)+d_{G'-T}(S)-r'(S),
\end{eqnarray*}
which implies that $G'$ has no fractional $r'$-factor. (Otherwise,
$r'(A)+d_{G'-A}(B)-r'(B)\geq0$ for all disjoint subsets $A$ and $B$ of $V(G)$ by Theorem 2.
Set $A=T$ and $B=S$.
Thus, we obtain $r'(T)+d_{G'-T}(S)-r'(S)\geq0$, a contradiction.) And so, $G$ has no fractional
$r'$-factor excluding $H$, that is, $G$ has no fractional $r$-factor including $H$. Hence, $G$
has no all fractional $(g,f)$-factors excluding $H$, a contradiction. This finishes the proof of
Theorem 4. \hfill $\Box$

\section{The proof of Theorem 5}
\noindent{\it Proof of Theorem 5.} \ According to Theorem 4, we need only to verify that
$$
g(S)+d_{G-S}(T)-f(T)\geq d_H(S)-e_H(S,T)
$$
for all disjoint subsets $S$ and $T$ of $V(G)$.

If $T=\emptyset$, then we have
$$
g(S)+d_{G-S}(T)-f(T)=g(S)\geq d_H(S)=d_H(S)-e_H(S,T).
$$
In the following, we assume that $T\neq\emptyset$. Note that $(g(x)-d_H(x))d_G(y)\geq(d_G(x)-d_H(x))f(y)$
holds for any $x,y\in V(G)$, that is, $g(x)d_G(y)\geq d_G(x)f(y)+d_H(x)(d_G(y)-f(y))$ holds for any
$x,y\in V(G)$. Hence, we have
$$
\Big(\sum_{x\in S}g(x)\Big)\Big(\sum_{y\in T}d_G(y)\Big)\geq\Big(\sum_{x\in S}d_G(x)\Big)\Big(\sum_{y\in T}f(y)\Big)+\Big(\sum_{x\in S}d_H(x)\Big)\Big(\sum_{y\in T}(d_G(y)-f(y))\Big),
$$
that is,
$$
g(S)d_G(T)\geq d_G(S)f(T)+d_H(S)(d_G(T)-f(T)).\eqno(2)
$$

We write $U=V(G)\setminus(S\cup T)$. Then we obtain
\begin{eqnarray*}
d_G(S)&=&e_G(S,T)+e_G(S,S)+e_G(S,U)\\
&\geq&e_G(S,T)+e_H(S,S)+e_G(S,U)\\
&=&e_G(S,T)+d_H(S)-e_H(S,T)-e_H(S,U)+e_G(S,U)\\
&\geq&e_G(S,T)+d_H(S)-e_H(S,T)\\
&=&d_G(T)-d_{G-S}(T)+d_H(S)-e_H(S,T),
\end{eqnarray*}
which implies
$$
d_G(S)-d_G(T)\geq-d_{G-S}(T)+d_H(S)-e_H(S,T).\eqno(3)
$$

In terms of (2) and (3), we have
\begin{eqnarray*}
&&d_G(T)(g(S)+d_{G-S}(T)-f(T)-d_H(S)+e_H(S,T))\\
&& \ \ \ =d_G(T)g(S)+d_G(T)d_{G-S}(T)-d_G(T)f(T)-d_G(T)d_H(S)+d_G(T)e_H(S,T)\\
&& \ \ \ \geq d_G(S)f(T)+d_H(S)(d_G(T)-f(T))+d_G(T)d_{G-S}(T)-d_G(T)f(T)\\
&& \ \ \ \ \ \ -d_G(T)d_H(S)+d_G(T)e_H(S,T)\\
&& \ \ \ =f(T)(d_G(S)-d_G(T))+d_G(T)d_{G-S}(T)-d_H(S)f(T)+d_G(T)e_H(S,T)\\
&& \ \ \ \geq f(T)(-d_{G-S}(T)+d_H(S)-e_H(S,T))+d_G(T)d_{G-S}(T)-d_H(S)f(T)+d_G(T)e_H(S,T)\\
&& \ \ \ =(d_{G-S}(T)+e_H(S,T))(d_G(T)-f(T))\geq0.
\end{eqnarray*}
Combining this with $d_G(T)\geq f(T)\geq|T|\geq1$, we obtain
$$
g(S)+d_{G-S}(T)-f(T)\geq d_H(S)-e_H(S,T).
$$
Theorem 5 is proved. \hfill $\Box$

\end{CJK*}
\end{document}